\documentclass{amsart}
 \newtheorem{thm}{Theorem}[section]

 \theoremstyle{definition}
 
 \theoremstyle{remark}

 \numberwithin{equation}{section}

\begin{document}

%
%
%
%
%
%
%
%
%

\title[Complexity of primes and squarefree integers]{Alternating state complexity of the set of primes and squarefree integers}
\author[J.-C Schlage-Puchta]{Jan-Christoph Schlage-Puchta}

\address{%
Ulmenstraße 69\\
18057 Rostock\\
Germany}

\email{jan-christoph.schlage-puchta@uni-rostock.de}

\subjclass{Primary 68Q45; Secondary 11N36}

\keywords{Automata, Complexity, Large Sieve}

\date{January 1, 2004}

\begin{abstract}
We show that the set of prime numbers has exponential alternating complexity, proving a conjecture by Fijalkow. We further show that the set of squarefree integers has essentially maximal possible alternating complexity.
\end{abstract}

\maketitle
\section{Introduction and results}
A deterministic automaton over an alphabet $\mathcal{A}$ is a set of states $S$, a starting position $s_0\in S$, a set of accepting states $A\subseteq S$, and a transition function $\delta:S\times\mathcal{A}\rightarrow S$. Given a word $w=a_1a_2\dots a_k\in\mathcal{A}^*$, we say that the word $w$ reaches the state $s_w=\delta(a_k, \delta(a_{k-1}, \dots, \delta(a_1, s_0)\dots)$, and we say that the automaton accepts $w$, if $s_w\in A$.

A non-deterministic automaton the function $\delta$ maps pairs $(a,s)$ to subsets of $Q$. Every sequence of choices of a state in $\delta(a,s)$ represents a valid run of the automaton, we say that a word is accepted, if there exists a run which ends in an accepting state. 

Chandra, Kozen and Stockmeyer \cite{CKS} defined an alternating automaton as an automaton, where the transition function maps pairs $(a,s)$ to finite positive boolean combinations of atomic formulae of the form $\delta=s'$ for some $s'\in S$. Here a positive boolean combination is a boolean function that can be written using $\wedge$ and $\vee$, but not $\neg$. A word is accepted, if Alice wins the following game: The game begins in $s_0$. When a letter $a$ is read, and the game is in state $s$, the boolean expression $\phi=\delta(a,s)$ is computed. If $\phi=(\delta=s')$ is atomic, the game continues to state $s$. If $\phi$ decomposes as $\psi\wedge\theta$, Eve replaces $\phi$ by one of $\psi, \theta$. If $\phi$ decomposes as $\psi\vee\theta$, Alice replaces $\phi$ by one of $\psi, \theta$. Alice wins, if after reading all fo $w$, the automaton reaches an accepting state, otherwise Eve wins. 

For example, if $\delta(a,s)=(\delta=s_1\wedge (\delta=s_2\vee\delta=s_3)$, the word $aw$ is accepted when starting in the state $s$ if and only if the word $w$ is accepted when starting in $s_1$, and the word $w$ is accepted when starting in at least one of the states $s_2$ and $s_3$.

Let $A$ be an infinite automaton. Karp \cite{Karp} defined the state complexity $f$ of $A$ as the function that associates to each integer $n$ the number of states of $A$ that are reachable within $n$ steps. We define the state complexity of a language $L$ as the function that associates to an integer $n$ the minimal number of states that an automaton recognizing $L$ can have, which are reachable by words of length $\leq n$. Note that the state complexity depends on the class of automata considered, so we have a deterministic, a non-deterministic and an alternating state complexity. Clearly, deterministic complexity is at least as large as non-deterministic complexity, which in turn bounds alternating complexity.

Fijalkow \cite{Alternating} showed that the alternating state complexity of the binary representations of the prime numbers is at least linear, and conjectured that it is exponential. Here we prove the following.

\begin{thm}
\label{thm:main}
We identify integers with words in $\{0,1\}^*$ by means of the binary expansion read from right to left.
\begin{enumerate}
\item The set of prime numbers has alternating complexity at least $c\frac{2^{n/2}\log n}{n}$.
\item The set of squarefree numbers has deterministic and alternating complexity $2^{n+\mathcal{O}(1)}$.
\end{enumerate}
\end{thm}

To give lower bounds for the automatic complexity we use an approach by Fijalkow \cite{Alternating}. Let $L$ be a language. We define the profile of a word $w$ of order $n$ as the set of all words $u$ of length $n$, such that $uw\in L$. We define the query table of $L$ of order $n$ as the set of all profiles of order $n$ of all words $w$. Define $q_L(n)$ as the size of the query table of $L$ of order $n$. This function is connected to the state complexity by means of the following.

\begin{thm}[Fijalkow]
\label{thm:query}
Let $L$ be a language of alternating state complexity $f$. Then $q_L(n)=2^{\mathcal{O}(f(n))}$.
\end{thm}
Hence, to find a lower bound for $f$ it suffices to find a lower bound for $q$.

\section{Primes}

For primes we use a counting argument to show that there have to be many different profiles. Fix an integer $n$, and put $N=2^n$, $x=\exp(\sqrt{N})$. For each integer $k$ consider the set
\[
Q_k=\left\{\mathcal{A}\subseteq[N]: \exists y<\frac{x}{N}\forall a<N: yN+a \mbox{ prime }\Leftrightarrow a\in\mathcal{A}\right\}.
\]
Clearly $Q_k$ is a subset of the query table of the set of primes of order $n$. By the prime number theorem the number of primes below $x$ is asymptotically equal to $\frac{x}{\log x}$. Call an interval $[yN, y(N+1)]$ rich, if it contains at least $\frac{N}{2\log x}$ primes, and poor otherwise. The poor intervals contain at most
\[
\frac{x}{N}\cdot\frac{N}{2\log x} =\frac{x}{2\log x}
\]
primes altogether. On the other hand a rich interval different from the interval $[0, N]$ contains at most as many primes as there are odd integers in this interval, hence, for $x$ sufficiently large the number of rich intervals is at least
\[
\left(\frac{1}{2}-\epsilon\right)\frac{x}{(N/2)\log x}-1 >\frac{x}{2N\log x}.
\]
Put $Q'=\bigcup_{k\geq \frac{N}{2\log x}} Q_k$. For each $\mathcal{A}\subseteq[N]$, let $\mathcal{Y}_{\mathcal{A}}$ be the set of all integers $y$ such that $\mathcal{A}=\{a<N: yN+a\mbox{ prime}\}$. Then the number of rich intervals equals 
\[
\sum_{k\geq \frac{N}{2\log x}} \sum_{\mathcal{A}\in Q_k} |\mathcal{Y}_{\mathcal{A}}|.
\]
Hence, 
\[
q_{\mbox{\footnotesize primes}}(n)\geq|Q'|\geq\frac{x}{2N\log x}\left(\max_{\vec{a}\subseteq\{1, \ldots, N\}} |\mathcal{Y}_a|\right)^{-1}
\]
We now bound $|\mathcal{Y}_a|$ from above using the following result by Elsholtz \cite{Christian}, which is based on the large sieve.
\begin{thm}[Elsholtz]
For all $c, \epsilon>0$ there exists some $x_0$, such that for all $x>x_0$ and all tuples $1\leq a_1<a_2<\dots<a_k<x$ with $k>c\log x$, the number of all integers $n\leq x$ such that $n+a_i$ is prime for all $i$ is bounded by 
\[
x\exp\left(-\left(\frac{1}{4}-\epsilon\right)\frac{\log x\log\log\log x}{\log\log x}\right).
\]
\end{thm}
We can apply this theorem to bound $|\mathcal{Y}_{\mathcal{A}}|$, since $\frac{N}{2\log x} = \frac{N}{2\sqrt{N}} =\frac{\log x}{2}$, and we obtain
\[
|\mathcal{Y}_a| \leq x\exp\left(-\left(\frac{1}{4}-\epsilon\right)\frac{\log x\log\log\log x}{\log\log x}\right).
\]
Hence, for $x>x_0(\epsilon)$ we obtain
\begin{eqnarray*}
Q' & \geq & \frac{x}{2N\log x} x^{-1}\exp\left(\left(\frac{1}{4}-\epsilon\right)\frac{\log x\log\log\log x}{\log\log x}\right)\\
 & > & \exp\left(\left(\frac{1}{4}-2\epsilon\right)\frac{\log x\log\log\log x}{\log\log x}\right).
\end{eqnarray*}
Inserting $\log x = \sqrt{N} = 2^{n/2}$ we obtain
\[
q_L(n) \geq \exp\left(\left(\frac{1}{4}-2\epsilon\right)\frac{2^{n/2}(\log n + \mathcal{O}(1))}{n\log 2/2}\right)>\exp\left(\frac{2^{n/2}\log n}{2n}\right).
\]
The first part of Theorem~\ref{thm:main} now follows from Theorem~\ref{thm:query}.

\section{Squarefree integers}

Note first that every set of integers can be encoded as a subset of nodes of the infinite binary tree, hence, every set of integers has deterministic complexity at most $2^{n+1}$. 

We now turn to the lower bound. Fix an integer $n$, and put $N=2^n$. Let $p_k$ be the largest prime number such that $p_k^2<N$, and put $q=\prod_{i=1}^k p_i^2$. Define $\mathcal{S}=\{a: a<N, a\mbox{ squarefree}\}$. We now construct for every subset $\mathcal{T}\subseteq\mathcal{S}$ an integer $y$ such that $yNq+a$ is squarefree if and only if $a\in\mathcal{T}$. 

Let $\ell$ be a sufficiently large parameter determined later, $q'=\prod_{i\leq \ell} p_i^2$. Define the residue class $y_0\pmod{q'}$ by putting $y_0N\equiv -a\pmod{p^2_{k+a}}$, if $a\in\mathcal{S}\setminus\mathcal{T}$, and $y_0\equiv 0\pmod{p_j^2}$ for all $j\leq \ell$ not of the form $k+a$ with $a\in\mathcal{S}\setminus\mathcal{T}$. We claim that if $\ell$ is sufficiently large, then there exists an integer $y$ satisfying $y\equiv y_0\pmod{q'}$ such that $yN+a$ is squarefree for $a\leq N$ if and only if $a\in\mathcal{T}$. 

If $a\not\in\mathcal{S}$, and for $y\equiv y_0\pmod{q'}$ we have $y\equiv 0\pmod{q}$, hence, $(y+a, q)=(a,q)$, and as $a$ is not squarefree, it is divisible by the square $p_i^2$ of a prime $<N$. We conclude that $p_i^2$ divides $y+a$, and $y+a$ is not squarefree. If $a\in\mathcal{S}\setminus\mathcal{T}$, then $y_0+a\equiv 0\pmod{p_{k+a}^2}$, hence, $y_0+a$ is not squarefree as well. We conclude that for all $y\equiv y_0\pmod{q'}$ and all $a\not\in\mathcal{T}$ we have that $y+a$ is not squarefree.

 We now show that for many $y$ we have that $a\in\mathcal{T}$ implies that $y+a$ is squarefree. Write $y=zq'+y_0$. If $i\leq k$, then $p_i^2|y$, therefore $p_i^2$ does not divide $y+a$ for any $a\in\mathcal{S}$. If $i>k$, then $p_i^2$ divides at most one integer in $[y, y+N]$. If $i\leq \ell$, then by construction of $y_0$ this integer is not in $\mathcal{T}$. If $i>\ell$, then there are at most $|\mathcal{T}|\left(\frac{Z}{p_i^2}+1\right)$ integers $z<Z$, such that there is some $a\in\mathcal{T}$ with $p_i^2|zq'+y_0+a$. Furthermore, if $p_i^2>(Z+1)q'+N$, then no such $z<Z$ exists. We conclude that the number of $z<Z$ such that $zq'+y_0+a$ is squarefree for all $z<Z$ is at least
 \[
Z -  \sum_{\ell<i<\sqrt{Zq'}} |\mathcal{T}|\left(\frac{Z}{q_i^2}+1\right) > Z - NZ\sum_{n=\ell+1}^\infty \frac{1}{n^2} - \sqrt{Zq'}.
 \]
 If we pick $\ell=2N$, this quantity is at least $\frac{Z}{2}-\sqrt{Zq'}$. We now choose $Z>4q'$, and find that the number of suitable $z$ is positive. We have constructed an integer with the desired profile,.
 
 We conclude that $q_L(n)\geq 2^{|\mathcal{S}|}$. It is well known that the set of squarefree numbers has density $\frac{6}{\pi^2}$, hence, 
 \[
 \log q_L(n)\geq |\mathcal{S}|\log 2 \geq \left(\frac{6\log 2}{\pi^2}-\epsilon\right)N = 2^{n+\mathcal{O}(1)}.
 \]
 Theorem~\ref{thm:query} now implies that the alternating complexity of the set of squarefree numbers is $\geq 2^{n-C}$ for some constant $C$. As the alternating complexity of a set is always bounded by the deterministic complexity, the second claim of Theorem~\ref{thm:main} follows.

\end{document}